\documentclass{article}

\usepackage{lineno,hyperref}
\usepackage{epsfig}
\usepackage{graphicx}
\usepackage{epstopdf}
\usepackage{amssymb,amsmath}
\usepackage{color}

\newtheorem{examp}{Example}

\newtheorem{proposition}{Proposition}
\newtheorem{corollary}{Corollary}

\def\C{\mathbb{ C}}

\newcommand{\R}{\mathbb{R}}

\begin{document}

\title{Filtering the Tau method with Frobenius-Pad\'e Approximants\footnote{This work was partially supported by CMUP(UID/MAT/00144/2013), which is funded by FCT (Portugal) with national and European structural funds (FEDER), under the partnership agreement PT2020}}

\author{Jo\~ao Carrilho de Matos\footnote{jem@isep.ipp.pt, Instituto Superior de Engenharia do Porto, Rua Dr. Ant\'onio Bernardino de Almeida, 431, 4249-015 Porto, Portugal}, Jos\'e M. A. Matos\footnote{Instituto Superior de Engenharia do Porto and Centro de Matem\'atica da Universidade do Porto}, Maria Jo\~ao Rodrigues\footnote{Faculdade de Ci\^encias da Universidade do Porto and  Centro de Matem\'atica da Universidade do Porto}}

\maketitle

\begin{abstract}
In this work, we use rational approximation to improve the accuracy of spectral solutions of differential equations. When working in the vicinity of solutions with singularities, spectral methods may fail their propagated spectral rate of convergence and even they may fail their convergence at all. We describe a Pad\'e approximation based method to improve the approximation in the Tau method solution of ordinary differential equations. This process is suitable to build rational approximations to solutions of differential problems when their exact solutions have singularities close to their domain. 
\end{abstract}

\textbf{keywords:}
Tau method, Pad\'e approximation, Froissart doublets.

\section{Introduction}
It is well known that spectral methods are very efficient to solve differential equations with smooth solutions without singularities close to the interval of orthogonality. In fact they exhibit exponential rate of convergence \cite{Canuto}.
However when the solution of a differential problem has singularities close to or on the interval of orthogonality spectral methods usually loose their efficiency. In fact, when the solution has singularities near to the orthogonality interval, the convergence of spectral methods is slow and when a solution has singularities on the interval of orthogonality spectral methods have only algebraic rate of convergence. There are several methods to improve the approximation given by the spectral solution, e.g. using extrapolation methods \cite{Brez,sidi2003}, filtering functions \cite{Canuto}, changing of variables or by decomposing the domain \cite{Peyret}. These post-processing methods are frequently called filtering processes  of  a spectral solution. 

In this paper we suggest to filter a tau solution  of a differential problem with a slow rate of convergence, using Chebyshev-Pad\'e or Legendre-Pad\'e approximants.
This choice is motivated by the theoretic results related to meromorphic functions \cite{Suetin82} and Markov type functions \cite{Gonchar92}.
These results are related to convergence acceleration and the domain expansion given by the partial orthogonal series. Moreover, it is possible to extract more information from Pad\'e approximants. In fact, we can determine singularities of functions using the poles
of Pad\'e approximants \cite{buslaev06,buslaev09}. Here we have to be aware of the fact that all the theoretic results mentioned above are related with Pad\'e approximants computed with the coefficients of a formal orthogonal expansion. In this filtering process we use the Tau coefficients which do not coincide with the orthogonal expansion coefficients. The Tau coefficients are affected by the errors inherent to the Tau method and by errors caused by the use of finite arithmetic. We will give a special relevance to the numeric errors since they origin Pad\'e approximants with Froissart doublets, which destroy the structure of the numeric Pad\'e table.

The purpose of this paper is to present some numerical results of application of this filtering process to several cases. In sections $2$ and $3$ we give the notations and algorithms concerning the tau method and Pad\'e approximation, respectively. In sections 4 we present the filtering process, in section  5 we obtain some properties of the filters  and we conclude in section 6 with comments and conclusions.

\section{The Operational Tau Method}

Here we describe an improved version \cite{JMatos2014} of the operational Tau method to solve linear ordinary differential equations  \cite{ortiz69,Ortiz81}.

Let $\mathcal{D}_\nu$ be the class of linear differential operators of order $\nu$ with polynomial coefficients and $D\in \mathcal{D}_\nu$
\begin{equation}\label{Ddef}
D \equiv \sum_{i=0}^\nu p_i (t) \frac{d^i}{dt^i},\ \ p_i(t)\in \mathbb{P}. 
\end{equation}
In order to solve a differential equation
\begin{equation} \label{Diffeqn}
D y(t) = f(t),\ t\in ]a,\ b[\subset\mathbb{R}
\end{equation}
Ortiz and Samara  \cite{Ortiz81} developed the operational approach for the Lanczos Tau method based on the algebraic representation of the linear differential operator \eqref{Ddef}. They proved that
\begin{equation*}
D\, y = \mathbf{y}\, \boldsymbol{\Pi}\, \mathbf{t}^T,\quad \text{with}\quad  \boldsymbol{\Pi}=\sum_{i=0}^\nu \boldsymbol{\eta}^i\, p_i (\boldsymbol{\mu}) 
\end{equation*}
where $y = \mathbf{y}\mathbf{t}^T=\sum_{i=0}^\infty y_i t^i$,  with $\mathbf{y}=[y_0,\ y_1,\ \ldots]$ and $\mathbf{t}=[1,\ t,\ t^2,\ \ldots]$, is the matricial representation of $y$, and the matrices $\boldsymbol{\eta}$ and $\boldsymbol{\mu}$ represent, respectively, the differentiation and the shift effects on $y$. 

Let $\{\phi_i\}_{i\geq 0}$ be a family of functions defined in an interval $I\subset\R$, orthogonal with respect to a weight function $w$, and let
\begin{equation*}\label{series1}
y = \sum_{i=0}^{\infty} c_{i}\phi_i,\quad c_{i}=\frac{\left\langle  y,\phi_i\right\rangle_w}{\left\|\phi_i\right\|_{w}^{2}},
\end{equation*}
be a formal series, where $\left\langle f,g \right\rangle_{w} =\int_{I}f(t)g(t)w(t)\text{d}t$ and $\left\| f\right\|_{w}=\sqrt{\left\langle f,f\right\rangle_{w}}$. 
 If $\{\phi_i\}_{i\geq 0}$ is a polynomial basis such that $\phi_i$ is a polynomial of degree $i$ and $\mathbf{V}=[v_{i,j}]_{i,j\geq 0}$ is the  matrix of the coefficients of those polynomials, that is, if 
\[ \phi_i = \sum_{j=0}^i v_{i j}t^j,\ i=0,1,\ldots \]
and if 
\begin{equation}\label{yv}
 y=\mathbf{y_v}\mathbf{v}^T=\sum_{i=0}^\infty c_i \phi_i,\quad \text{with}\quad \mathbf{v}^T=\mathbf{V}\mathbf{t}^T. 
\end{equation} 
then the effect of the differential operator \eqref{Ddef} in the coefficients of $y$ can be represented by the effect of an algebraic operator over the vector $\mathbf{y_v}$, and is given by \cite{Ortiz81}
$$ 
D y=\mathbf{y_v}\, \boldsymbol{\Pi_{\phi}} \mathbf{v}^T,\quad \text{where}\quad \boldsymbol{\Pi_{\phi}} = \mathbf{V}\, \Pi\, \mathbf{V}^{-1} 
$$

\subsection{Improved operational Tau method}

At this point, we must make two remarks: the first one is that, in practice, numerical methods work with finite matrices and the second one is that
this procedure can be numerically unstable, particularly if the condition number of $\mathbf{V}$ is large.
However, if $\{\phi_i\}_{i\geq 0}$ is an orthogonal polynomials basis,  Matos et al presented in \cite{JMatos2014} a  method to overcome this drawback. 
Since  the polynomials $\phi_{i}$  satisfy a three term recurrence relation

\begin{equation}\label{RecRel}
\left\{\begin{array}{ll}
t \phi_i(t)=\alpha_i \phi_{i+1}(t)+\beta_i \phi_i(t)+\gamma_i \phi_{i-1}(t),  & i=1,2,\ldots \\
\phi_0(t)=1,\ \phi_1(t)=(t-\beta_0)/\alpha_0
\end{array}\right.
\end{equation}
then the matricial representation of the shift effect takes the form
$$
t\, y =\mathbf{y_v}\, \boldsymbol{\mu_{\phi}}\, \mathbf{v}^T
$$
where
\begin{equation*}
\boldsymbol{\mu_{\phi}}=\left[\begin{array}{cccc} \beta_0 & \alpha_0 \\ \gamma_1 & \beta_1 & \alpha_1 \\ & \gamma_2 & \beta_2 & \alpha_2 \\ & & & \cdots \end{array}\right]
\end{equation*}
On the other hand, defining $\eta_{i,j}$ as the coefficients in
\begin{equation*}\label{DifPi}
\frac{d}{dt}\phi_i = \sum_{j=0}^{i-1}\eta_{i j}\phi_j,\ i= 0,1,\ldots
\end{equation*}
and defining $\eta_\phi =[\eta_{i j}]_{i,j\geq 0}$, then
\[
\frac{d}{dt} y =\mathbf{y_{v}}\, \mathbf{\eta_{\phi}}\, \mathbf{v}^T.
\]
\noindent By differentiating both sides of recurrence relation \eqref{RecRel}, it is easy to  see that the coefficients $\eta_{i j}$
satisfy the following recurrence relation

\begin{equation*}
\left\{\begin{array}{l}
\eta_{i+1,j}=\frac{1}{\alpha_i}(\alpha_{j-1}\eta_{i,j-1}+(\beta_j-\beta_i)\eta_{i,j} +\gamma_{j+1}\eta_{i,j+1} - \gamma_i \eta_{i-1,j}),\ j=0,\ldots,i-1 \\
\eta_{i+1,i}=\frac{1}{\alpha_i}(\alpha_{i-1}\eta_{i,i-1}+1)
\end{array}\right.
\end{equation*}
\noindent for $i=1,2,\ldots$ and
\[
\left\{\begin{array}{ll}
\eta_{0,j}=0, & j=0,1,\ldots \\
\eta_{1,0}=\frac{1}{\alpha_0}, \\
\eta_{1,j}=0, & j=1,2,\ldots
\end{array}\right.
\] 

Having defined $\boldsymbol{\mu_{\phi}}$ and $\boldsymbol{\eta_{\phi}}$, we get 

\begin{equation}\label{piphi}
D\, y = \mathbf{y_v}\, \boldsymbol{\Pi_{\phi}}\, \mathbf{v}^T,\ \boldsymbol{\Pi_{\phi}} = \sum_{i=0}^\nu \boldsymbol{\eta}^{i}_{\phi}\, p_i (\boldsymbol{\mu}_{\phi}) 
\end{equation}
Thus, we don't need to compute the inverse of $\mathbf{V}$, which  stabilizes the operational Tau method.

The operational approach of the Tau method  is based on the matrix $\Gamma_{\phi}=[G\quad \overline{\Pi}_{\phi}]$ where, for some $n\in\mathbb{N},\ n\geq \nu$, $G\in \mathbb{R}^{(n+1)\times \nu}$ is the matrix representation of the initial conditions of the differential problem and $\overline{\Pi}_{\phi}\in \mathbb{R}^{(n+1)\times (n+1-\nu)}$ is the matrix $\Pi_{\phi}$ truncated to its first lines and first columns \cite{JMatos2014, Ortiz81}. If $\Gamma_{\phi}$ is a regular matrix then we solve an algebraic system of linear equations in order to obtain $y_n=[c_0^{(n)},\ldots,c_n^{(n)}]$, the coefficients of the Tau approximant of $y_v$ in (\ref{yv}).

\subsection{Numerical example}
In order to explain the proposed filtering method, we begin by considering the following example, where some useful notation is introduced.

\begin{examp}\label{exampleSQRT}

Here, we will consider the function $y(t)=\frac{\pi}{4}\sqrt{2(t+1)}$, which has a branch cut on $]-\infty,-1]$. It is known that the Chebyshev expansion of this function \cite{Paszkowski84}
is
\begin{equation}
y = \sum_{k=0}^{\infty}c_{k}T_{k} = 1+\sum_{k=1}^{\infty}(-1)^{k+1}\frac{2}{4k^{2}-1}T_{k}. \label{yex1}
\end{equation}

The function $y$ is analytic on $\C\setminus ]-\infty , -1]$ thus, the singularity, $\zeta$, closest to the interval of orthogonality, $[-1,1]$, coincides with the extreme point $\zeta=-1$.

We can define $y$ as the solution of the linear ordinary differential equation
\begin{equation}\label{ex1tau}
	(t+1)\frac{\text{d}y}{\text{d}t}-\frac{1}{2}y=0
\end{equation}
\noindent with condition $y(0)=\frac{\pi}{4}\sqrt{2}$.
\end{examp}

From (\ref{piphi}), for the differential operator $D$ associated to equation (\ref{ex1tau}), we get $\Pi_{\phi}=\eta_{\phi}(\mu_{\phi}+I)-\frac{1}{2}I$, where $I$ is the infinite identity matrix. The initial condition is translated into the matrix $G=(\phi_0(0),\phi_1(0),\ldots)^T$. Since we are considering $\phi_i=T_i,\ i\geq 0$, the Chebyshev polynomials, then
\[
\Gamma_T= \left[
\begin{array}{rccccc}
1 & -1/2 & \\
0 & 1 & 1/2 \\
-1 & 2 & 4 & 3/2\\
0 & 3 & 6 & 6 & 5/2 \\
1 & 4 & 8 & 8 & 8 & 7/2 \\
&&&&&\cdots
\end{array}
\right]
\]

\begin{figure}[hbt]
\begin{center}
\includegraphics[width=1\textwidth]{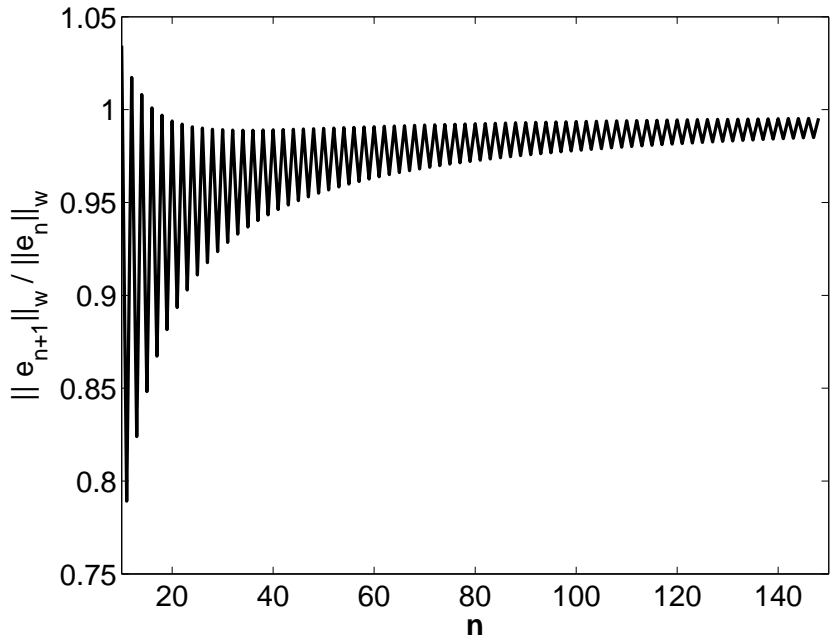}
\caption{  Rate of consecutive errors $\| e_{n+1}\|_{w}/ \|e_{n}\|_{w}$, for $n=9,10,\ldots,149$.}
\label{fig1tauconv}
 \end{center}
\end{figure}

As the singularity $\zeta$ coincides with the extreme point of the interval of orthogonality, the convergence of the Tau method is slow. To emphasize this fact, we show in Figure \ref{fig1tauconv}, the rate of consecutive errors $\| e_{n+1}\|_{w}/ \|e_{n}\|_{w}$, for $n=9,10,\ldots,149$, where $e_{n}=y-y_{n}$, and $\left\| e_{n} \right\|_{w}^{2} = \int_{-1}^1 e_{n}^{2}(t)w(t) \text{d}t$ with $w(t)=(1-t^{2})^{-1/2}$. 
We can observe that  the rate $\| e_{n+1}\|_{w}/ \|e_{n}\|_{w}$ approaches 1 when $n$ increases. 

\subsection{Error analysis in the Tau Method}

As noted by \cite{Ortiz81} for the operational approach of the Tau method, the polynomial solution $y_n$, obtained by solving the linear algebraic system with the truncated matrix $\Gamma_{\Phi}$, is the exact solution of the perturbed differential equation
\begin{equation} \label{TauDiffeqn}
D y_n(t) = f(t) + \tau_n(t),\ t\in ]a,\ b[\subset\mathbb{R}
\end{equation}
where $\tau_n$ is the polynomial residual resulting from the truncation process of matrix $\Gamma_{\Phi}$. Since we are considering linear differential operators $D$, then, subtracting side by side (\ref{TauDiffeqn}) from (\ref{Diffeqn}) results that the error $e_n=y-y_n$ in the Tau method, is the exact solution of the differential equation
\begin{equation} \label{ErrorDiffeqn}
D e_n(t) = -\tau_n(t),\ t\in  ]a,\ b[\subset\mathbb{R}
\end{equation}
with homogeneous conditions.

\begin{figure}[hbt]	
\begin{center}
\includegraphics[width=1.1\textwidth]{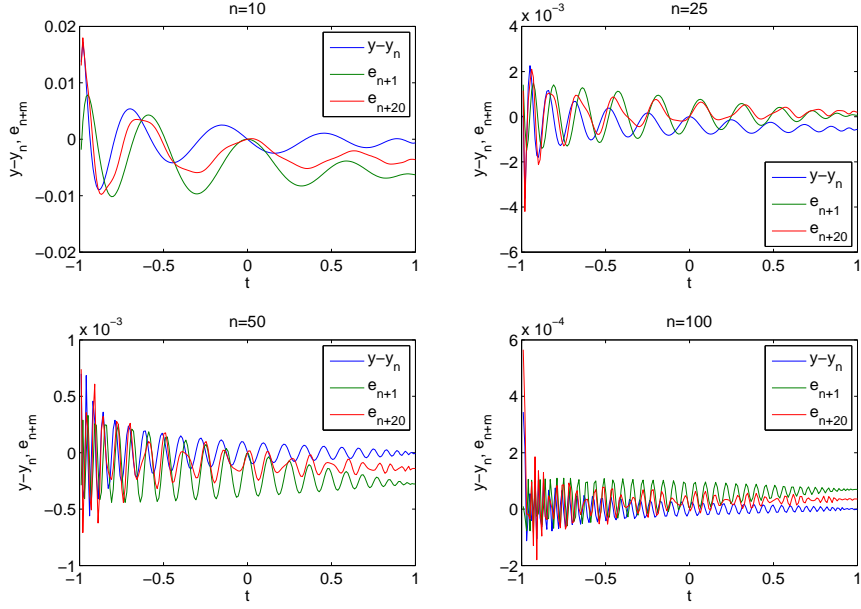}
\caption{Error curves $y(t)-y_n(t)$ and error curves approximations $\tilde{e}_{n+1}(t)$ and $\tilde{e}_{n+20}(t)$ obtained with the Tau method for example \ref{exampleSQRT}.}
\label{ErrorsAprox}
 \end{center}
\end{figure}

Based on that property, some authors \cite{MJRodrigues13} developed an \textit{a posteriori} error analysis, solving by the Tau method the  differential equation in (\ref{ErrorDiffeqn}) and getting a polynomial approximation of, say, degree $n+m$, as the exact solution of
\begin{equation} \label{TauErrorDiffeqn}
D \tilde{e}_{n+m}(t) = -\tau_n(t)+\tau_m(t),\ t\ ]a,\ b[\subset\mathbb{R}
\end{equation}
In Figure \ref{ErrorsAprox} we show, for a selected set of $n$ values, and with $m=1$ and $m=20$, the error curves $y(t)-y_n(t)$ and the curves $\tilde{e}_{n+1}(t)$ and $\tilde{e}_{n+20}(t)$. We can see the effective numerical estimation of the error, even for $m=1$, and that, in general, the numerical estimator $e_{n+20}$ follows the error closer than the estimator $e_{n+1}$.

For the filtering process, proposed in section \ref{FilterinProcess}, the coefficients errors in the Tau method
\[
\Delta c^{(n)}=\left( \Delta c^{(n)}_{0}, \Delta c^{(n)}_{1},\ldots ,\Delta c^{(n)}_{n}\right), \quad \Delta c^{(n)}_{k}= c_{k}- c^{(n)}_{k}, \ k=0,1,\ldots,n
\]
play a main role in the final results. We expect that, for each $k=0,1,\ldots,n$, the error $\Delta c^{(n)}_{k}$ decreases with increasing $n$, provided that the Tau method converges, and also that this errors do not affect too much the filtered solution error. In general cases, with the Tau method, those values $\Delta c^{(n)}_{k}$ are the coefficients of the error function $e_n=y-y_n$ and can be approximated by the coefficients of the a posteriori error estimator $\tilde{e}_{n+m}$.

\subsection{Errors on the coefficients in example \ref{exampleSQRT}}
In the previous example we can simplify and solve exactly the system associated to the Tau method, obtaining exact values for $c^{(n)}_k$. Writing the last $n-1$ equations in the form
\[ 
\sum_{i=k}^n 2i c^{(n)}_i + (k-\frac{3}{2})c^{(n)}_{k-1} =0,\ k=2,\ldots,n 
\]
and subtracting each equation $k-1$ from equation $k$, we get the backward recurrence
\begin{equation}\label{cnk}
c^{(n)}_k = -\frac{2k+3}{2k-1}c^{(n)}_{k+1},\ k=1,\ldots,n-2,\quad\text{and}\quad c^{(n)}_0 = \frac{3}{2}c^{(n)}_{1} 
\end{equation}
This result leads, by mathematical induction, to
\begin{equation}\label{cnk}
	\left\{ \begin{array}{ll}
	c^{(n)}_k = (-1)^{n-k}\frac{4n(2n-1)}{4k^2-1}c^{(n)}_{n}, & k=1,\ldots,n-1,\\ 
	c^{(n)}_0 = (-1)^{n+1} 2n(2n-1)c^{(n)}_{n} &
\end{array}\right.
\end{equation}
and, substituting in the first equation, we can solve for $c^{(n)}_{n}$, getting
\[
  c^{(n)}_n = (-1)^{n+1} \frac{y_0}{2n(2n-1)S_n} 
\]
where $y_0=y(0)$ is the initial value, and 
\[ 
S_n=1+\sum_{k=1}^{n} c_{k}T_{k}(0)-\frac{1}{2n(2n+1)}T_{n}(0) 
\] 
is the partial sum for $y(0)$ if $n$ is odd and is the partial sum plus a correction term in the last coefficient if $n$ is even. Substituting in (\ref{cnk}) we get
\[  
c^{(n)}_k = (-1)^{k+1} \frac{2y_0}{(4k^2-1)S_n},\ k = 1,\ldots,n-1,\quad\text{and}\quad c^{(n)}_0 = \frac{y_0}{S_n} 
\]
and so, comparing these coefficients  $c^{(n)}_k $ with the exact coefficients $c_k$ \eqref{yex1}, we verify that 
\[
c_n^{(n)}=(\frac{2n+1}{4n})\frac{y_0}{S_n}c_n\quad\text{and}\quad c_k^{(n)}=\frac{y_0}{S_n}c_k,\ k = 0,\ldots,n-1
\]
This means that, for fixed $n\in\mathbb{N}$, each Tau coefficient $c_k^{(n)}$, except the last one, is the exact coefficient $c_k$ times a constant factor. This is relevant, in next sections, for our filtering results and to justify the exact formula for the error in the Tau coefficients
\begin{equation}\label{Deltackn}
\Delta c_k^{(n)} = \left\{
\begin{array}{ll}
(1-\frac{y_0}{S_n})c_k, & k=0,\ldots,n-1 \\
(1-(\frac{2n+1}{4n})\frac{y_0}{S_n})c_n & k=n
\end{array}\right. 
\end{equation}

Another property, resulting from the last column values of the $\Gamma_{\phi}$ matrix, for this particular example, is that the residual is $\tau_n(t)=(n-\frac{1}{2})c_n^{(n)} T_{n+1}(t)$. Using the previous formula for $c_n^{(n)}$, we get 
\[ 
\tau_n(t)=(-1)^{n+1}\frac{y_0}{4n S_n} T_{n+1}(t),
\]
and so, the residual $\tau_n$ is approximating zero in $[-1,\ 1]$, with uniform norm and with amplitudes decreasing with $n$.

Before we proceed with the proposed filtering method, we will remember the definition of Frobenius-Pad\'e approximants, also known as linear Pad\'e approximants from series of orthogonal polynomials.

\section{Frobenius-Pad\'e approximants from orthogonal series}

We begin this section by defining Frobenius-Pad\'e approximants from orthogonal series.  
Let $\{\phi_i\}_{i\geq 0}$ be an orthogonal polynomial basis.
Given two nonnegative integers $p$ and $q$, we say that the rational function
\begin{equation}\label{phi}
\Phi_{p,q}(y)=\frac{N_{p,q}}{D_{p,q}}=\frac{\sum_{i=0}^{p}a_{i}\phi_i}{\sum_{i=0}^{q} b_{i}\phi_i}
\end{equation}

\noindent is a Frobenius-Pad\'{e} approximant of type $(p,q)$ from the  series $y$ \cite{ACMatos01} if
\begin{equation}\label{lindef}
D_{p,q} y - N_{p,q} = \sum_{i=p+q+1}^{\infty}e_{i}\phi_i.
\end{equation}

In order to determine the coefficients $a_i,\ i=0,1,\ldots,p$ and $b_i,\ i=0,1,\ldots,q$ in (\ref{phi}) we introduce $h_{j,i}, j=0,1,\ldots$ such that
\[
\phi_{i} y = \sum_{j=0}^{\infty}h_{j,i}\phi_{j},\quad i=0,1,\ldots.
\]
\noindent Thus $h_{j,i}$, j=0,1,\ldots, are the coefficients of the orthogonal series  $\phi_{i} y,\quad i=0,1,\ldots$.

Identifying coefficients in (\ref{lindef}) we can show that the coefficients $a_{j}$ and $b_{i}$ of $\Phi_{p,q}$ are solutions of the following homogeneous system of $p+q+1$ linear equations and $p+q+2$ unknowns
\begin{equation}\label{eq:afp1D11}
\left\{\begin{array}{ll}
	\displaystyle{\sum_{i=0}^{q}h_{j,i}b_{i}-a_{j}}=0, &  j=0,\ldots,p \\
	\displaystyle{\sum_{i=0}^{q}h_{j,i}b_{i}}=0, & j=p+1,\ldots,p+q. 
	\end{array}\right.
\end{equation}

\noindent which always admits a non trivial solution.

If we set $b_{q}=1$ then we can use equations \eqref{eq:afp1D11} to determine
the coefficients of the normalized approximant, based in the matricial form introduced in the following proposition \cite{ACMatos01}

\begin{proposition}\label{Prop:NormalizedApprox}
Let $p,q\in \mathbb{N}_0$,
\[ \mathbf{g}^{[p/q]}=\left[h_{0,q}\ldots h_{p,q} \right]^{T},\quad 
\mathbf{h}^{[p/q]}=\left[h_{p+1,q}\ldots h_{p+q,q} \right]^{T}, \]
and
\[ \mathbf{G}^{[p/q]}=\left[\begin{array}{ccc}
h_{0,0}&\cdots & h_{0,q-1}\\
\vdots & &\vdots \\
h_{p,0}&\cdots &h_{p,q-1}
\end{array} \right], \quad \text{and} \quad 
\mathbf{H}^{[p/q]}=\left[\begin{array}{ccc}
h_{p+1,0}&\cdots & h_{p+1,q-1}\\
\vdots & &\vdots \\
h_{p+q,0}&\cdots &h_{p+q,q-1}
\end{array} \right]. \]
If $\mathbf{H}^{[p/q]}$ is nonsingular then  
\[ \mathbf{a}=\left[ a_{0}\ldots a_{p} \right]^{T}, \quad 
\mathbf{b}=\left[ b_{0}\ldots b_{q-1} \right]^{T}, \]
are determined by

\begin{align}
& \mathbf{H}^{[p/q]}\cdot\mathbf{b}=-\mathbf{h}^{[p/q]}\label{AFPsist1}\\ 
& \mathbf{a}=\mathbf{G}^{[p/q]}\cdot\mathbf{b}+\mathbf{g}^{[p/q]}\label{AFPsist2}
\end{align}

\end{proposition}

In the conditions of this proposition, the coefficients of the denominator, $b_{i}$, $i=0,1,\ldots,q-1$,  are uniquely determined by solving \eqref{AFPsist1}. Once determined the coefficients $b_{i}$, we use \eqref{AFPsist2} to compute the numerators coefficients $a_{i}$, $i=0,1,\ldots,p$.

The characteristic recurrence relation  \eqref{RecRel}  of the orthogonal polynomials $\phi_{i}$ leads to the following proposition,  allowing the computation of the entries in $H^{[p/q]},\ G^{[p/q]},\ h^{[p/q]}$ and $g^{[p/q]}$

\begin{proposition}\label{Prop:ReqRel}
The coefficients $h_{i j}$ can be computed \cite{ACMatos01} using the recurrence relation

\begin{equation}\label{AFP:eqrel1}
	h_{i,j+1}=\frac{1}{\alpha_{j}}\left( \frac{\mu_{i+1}}{\mu_{i}}\alpha_{i}h_{i+1,j}+(\beta_{i}-\beta_{j})h_{i,j}+
	\frac{\mu_{i-1}}{\mu_{i}}\gamma_{i}h_{i-1,j}-\gamma_{j}h_{i,j-1}\right), \ i,j=1,2,\ldots
\end{equation}
\noindent and,
\[
h_{i,0}=c_{i}, \ \ i=0,1,\ldots, \ \ \ 
h_{0,j}=\frac{\mu_{j}}{\mu_{0}}h_{j,0},  \ \ \ \ j=1,2,\ldots, 
\]
where $h_{i,-1}=0$, $\alpha_{i},\beta_{i}$ and $\gamma_{i}$ are the coefficients in \eqref{RecRel} and $\mu_{i}=\left\|\phi_{i}(t)\right\|_{w}^{2}$.
\end{proposition}

In particular, for Chebyshev and for Legendre polynomials we have

\paragraph{\textbf{Chebyshev Polynomials}:} The Chebyshev polynomials, normalized with $T_{i}(1)=1, \ \ i=0,1,\ldots$, satisfy the recurrence relation \eqref{RecRel} with
\[
\left\{
\begin{array}{lll}
\alpha_{i}=\gamma_{i}=\frac{1}{2}, & \beta_{i}=0, & i=1,2,\ldots \\
\alpha_{0}=1, & \beta_{0}=0
\end{array}
\right.
\]
and $\mu_{0}=\pi$, $\mu_{i}=\pi /2,\ i=0,1,\ldots$. Thus, the recurrence relation \eqref{AFP:eqrel1} takes the form
$$
\left\{
\begin{array}{ll}
h_{i,0}=c_{i}, & i=0,1,\ldots \\
h_{0,j}=\displaystyle{\frac{1}{2}c_{j}}, & j=1,2,\ldots \\
h_{1,1}=\displaystyle{h_{0,0}+\frac{1}{2}h_{2,0}} \\
h_{i,1}=\displaystyle{\frac{1}{2}(h_{i-1,0}+h_{i+1,0})}, & i=2,3,\ldots \\
h_{1,j}=2h_{0,j-1}+h_{2,j-1}-h_{1,j-2}, & j= 2,3,\ldots \\ 
h_{i,j}=h_{i-1,j-1}+h_{i+1,j-1}-h_{i,j-2}, & i,j= 2,3,\ldots \\ 
\end{array}
\right.
$$

With these formulas we can build direct formulas for some sequences of Chebyshev-Pad\'e approximants.

\begin{corollary}\label{ChebyPadeApprox}
Let $ y = \sum_{k=0}^\infty c_k T_k $
be a formal Chebyshev series and let $\Phi_{p,q}(y)=N_{p,q}/D_{p,q}$, with
\[
D_{p,q} =T_q + \sum_{i=0}^{q-1}b_{i}T_i, \quad \text{and}\quad N_{p,q} = \sum_{i=0}^{p}a_{i}T_i
\]
be its $(p,q)$ Chebyshev-Pad\'e approximant, then 
\begin{itemize}
	\item[(a)] $\forall p\in\mathbb{N}_0$ such that $c_{p+1}\neq 0$ we have
	\[ D_{p,1}(t)=b_{0}+t,\quad\text{and}\quad N_{p,1}(t)=\frac{1}{2}\sum_{k=0}^p (c'_{k-1}+c_{k+1}+2b_{0}c_k)T_k(t) \]
	with 
	\[ b_{0}=-\frac{c'_p+c_{p+2}}{2c_{p+1}} \]
	where $c'_{-1}=0,\ c'_0=2c_0$ and $c'_k=c_k,\ k\geq 1$. 
	\item[(b)] $\forall p\in\mathbb{N}_0$ such that the determinant
\[ \underline{\Delta} = \left|
\begin{array}{cc}
	c_{p+1} & c'_{p}+c_{p+2} \\
	c_{p+2} & c_{p+1}+c_{p+3}	
\end{array} \right| \neq 0 \]	 
we have
	\[ D_{p,2}(t)=b_0+b_1 t+T_2(t),\quad\text{and}\quad N_{p,2}(t)=\sum_{k=0}^p a_k T_k(t) \]
	with 
	\[ a_k=c'_k b_0+\frac{1}{2}(c'_{k-1}+c_{k+1})b_1+\frac{1}{2}(c'_{k-2}+c_{k+2}),\ k=0,\ldots,p \] 
$$
b_{0}=-\frac{\left|
\begin{array}{cc}
c_{p-1}+c_{p+3} & c_{p}+c_{p+2}\\
c'_{p}+c_{p+4} & c_{p+1}+c_{p+3}	
\end{array}\right|}{2\underline{\Delta}}
 \quad \text{and} \quad
b_{1}=-\frac{\left|
\begin{array}{cc}
c_{p+1} & c_{p-1}+c_{p+3}\\
c_{p+2} & c'_{p}+c_{p+4}	
\end{array}\right|}{\underline{\Delta}}.
$$
	where $c'_{-2}=0,\ c'_{-1}=2c_1-c_0,\ c'_0=2c_0$ and $c'_k=c_k,\ k\geq 1$. 
\end{itemize}
\end{corollary}

\paragraph{\textbf{Legendre Polynomials}:} The Legendre polynomials, normalized with $P_{i}(1)=1, \ \ i=0,1,\ldots$, satisfy the recurrence relation \eqref{RecRel} with
\[
\alpha_{i}=\frac{i+1}{2i+1}, \quad \beta_{i}=0, \quad \gamma_{i}=\frac{i}{2i+1}, \quad
\text{and}\quad \mu_{i}=\frac{2}{2i+1},\quad i=0,1,\ldots 
\] 

Thus, \eqref{AFP:eqrel1} takes the form
$$
\left\{
\begin{array}{ll}
h_{i,0}=c_{i}, \qquad h_{0,i}=\displaystyle{\frac{1}{2i+1}c_{i}}, & i=0,1,\ldots \\
h_{i,1}=\displaystyle{\frac{i+1}{2i+3}c_{i+1}+\frac{i}{2i-1}c_{i-1}}, & i=1,2,\ldots \\
h_{i,j+1}=\displaystyle{\frac{2j+1}{j+1}\left[\frac{i+1}{2i+3}h_{i+1,j}+\frac{i}{2i-1}h_{i-1,j}\right]-\frac{j}{j+1}h_{i,j-1}}, & i,j=1,2,\ldots
\end{array}
\right.
$$

In that case, even not so simple as in the Chebyshev case, formulas for Legendre-Pad\'e approximants $\Phi_{p,1}$ and $\Phi_{p,2}$, in terms of the series coefficients $c_k$ can be derived.

\begin{corollary}\label{LegPadeApprox}
Let $ y = \sum_{k=0}^\infty c_k P_k $
be a formal Legendre series and let $\Phi_{p,q}(y_{n})=N_{p,q}/D_{p,q}$, with
\[
D_{p,q} = P_q + \sum_{i=0}^{q-1}b_{i}P_i, \quad \text{and}\quad N_{p,q} = \sum_{i=0}^{p}a_{i}P_i
\]
be its $(p,q)$ Legendre-Pad\'e approximant, then
\begin{itemize}
	\item[(a)] $\forall p\in\mathbb{N}_0$ such that $c_{p+1}\neq 0$ we have
	\[ D_{p,1}(t)=b_{0}+t,\quad\text{and}\quad N_{p,1}(t)=\sum_{k=0}^p \left(h_{k,1}+b_{0}c_k\right)P_k(t) \]
	with 
	\[ b_{0}=-\frac{h_{p+1,1}}{c_{p+1}},\quad\text{and}\quad h_{k,1}=\frac{k}{2k-1}c_{k-1}+\frac{k+1}{2k+3}c_{k+1},\ k=0,1,\ldots \]
	\item[(b)] $\forall p\in\mathbb{N}_0$ such that the determinant
\[ \underline{\Delta} = \left|
\begin{array}{cc}
	c_{p+1} & h_{p+1,1} \\
	c_{p+2} & h_{p+2,1}	
\end{array} \right| \neq 0 \]	 
we have
	\[ N_{p,2}(t)=b_0+b_1 t+P_2(t),\quad\text{and}\quad D_{p,2}(t)=\sum_{k=0}^p a_k P_k(t) \]
	with 
	\[ a_k=c_k b_0+h_{k,1}b_1+h_{k,2},\quad k=0,\ldots,p \] 
$$
b_{0}=-\frac{\left|
\begin{array}{cc}
	h_{p+1,2} & h_{p+1,1} \\
	h_{p+2,2} & h_{p+2,1}	
\end{array}\right|}{\underline{\Delta}}
 \quad \text{and} \quad
b_{1}=-\frac{\left|
\begin{array}{cc}
	c_{p+1} & h_{p+1,2} \\
	c_{p+2} & h_{p+2,2}	
\end{array}\right|}{\underline{\Delta}}.
$$
\noindent	where 

	\begin{multline*}
	h_{k,2}=\frac{3(k^2-1)}{2(2k+3)(2k-1)}\left[\left(\frac{4k+3}{(2k-3)(k+1)}+1\right)c_{k-2}\right. \\ \left.+\frac{2k}{3(k-1)}c_{k}+\left(\frac{3}{(2k+5)(k-1)}+1\right)c_{k+2}\right],\ k=0,1,\ldots
	\end{multline*}
\end{itemize}
\end{corollary}

In the next section we will describe some numerical problems related with this filtering process and we will use the example \ref{exampleSQRT} to illustrate them.

\section{The Filtering Process}\label{FilterinProcess}

Let $y$ be the solution of a given differential problem, with formal orthogonal expansion $y=\sum_{i\geq 0}c_{i}\phi_{i}$  and let
\[
\Phi_{p,q}(y)=\frac{N_{p,q}}{D_{p,q}}=\frac{\displaystyle{\sum_{i=0}^{p}a_{i}\phi_i}}{\displaystyle{\sum_{i=0}^{q-1}b_{i}\phi_i + \phi_q}}
\]
be its Pad\'e approximant of type $(p,q)$. Since we have not access to the exact coefficients $c_i$, our filtering process makes use of the coefficients $c_i^{(n)}$ of the Tau solution of order $n$, $y_n = \sum_{i=0}^{n}c^{(n)}_{i}\phi_{i}$, of the given differential problem to construct the Pad\'{e} approximant of type $(p,q)$
\[
\Phi_{p,q}(y_{n})=\frac{N^{(n)}_{p,q}}{D^{(n)}_{p,q}}=\frac{\displaystyle{\sum_{i=0}^{p}a^{(n)}_{i}\phi_i}}{\displaystyle{\sum_{i=0}^{q-1}b^{(n)}_{i}\phi_i + \phi_q}}, \quad  \text{with} \quad p+2q+1\leq n.
\]

\subsection{Error in the filtering process}

A filter $\Phi_{p,q}(y_{n})$ represents the exact solution $y$ with an error that depends on the  errors on the coefficients $\Delta c^{(n)}$ and also on numerical errors caused by the numerical instability of the algorithm described above to compute Pad\'e approximants. In fact, the critical step is the resolution of the system of linear equations \eqref{AFPsist1}. To be more precise, the matrices $\mathbf{H}^{\left[ p/q\right]}$ are ill-conditioned for sufficiently large values of $p$ and $q$.
We remark that, in example \ref{exampleSQRT},  both matrices $\mathbf{H}^{\left[ p/q\right]}$, built with Tau coefficients or built with Fourier series coefficients, have condition numbers with same order of magnitude.

These numerical errors, caused by the numerical instability of the algorithm, are a serious drawback in the filtering method.
In effect, they origin Pad\'e approximants with Froissart doublets located nearby of the real segment $I=[-1,1]$. This behaviour is similar
to the Chebyshev and Legendre Pad\'e approximants computed with expansions perturbed with random noises \cite{JCMatos2014,JCMatos2015}.
The presence of Froissart doublets apart from destroying the structure of the computed Pad\'e table also spoil locally the approximation
given by the Pad\'e approximants.
In order to bypass this drawback we build a table, that we call \textit{Froissart table}.

\subsection{The Froissart table}
The formal definition of Froissart doublet, \cite{Stahl98}, is given in an asymptotic way, being thus useless for our purposes.
 For practical effects we will consider a Froissart doublet of a Pad\'e approximant as a pair $(\zeta,\eta)$ such that $\zeta$ is a pole, $\eta$ is a zero and $|\zeta-\eta|<\text{\rm tol}$, where $\text{\rm tol}>0$ is a prescribed \textit{tolerance}.
The $(p,q)$ entry of the $\text{\rm tol}$-Froissart table, $n_{p,q}$, is found by computing the zeros and poles of $\Phi_{p,q}$
and $n_{p,q}$ is the number of pairs of poles/zeros at distance less than $\text{\rm tol}$, or, by other words $n_{p,q}$ is the number
of Froissart doublets of $\Phi_{p,q}$.
We note that numerically $\Phi_{p,q}$ is a rational function with numerator of \textit{numeric} degree $p-n_{p,q}$
and denominator of \textit{numeric} degree $q-n_{p,q}$, since we have $n_{p,q}$ pairs of factors that \textit{almost} cancel.
\begin{figure}[hbt]
\begin{center}
\includegraphics[width=1\textwidth]{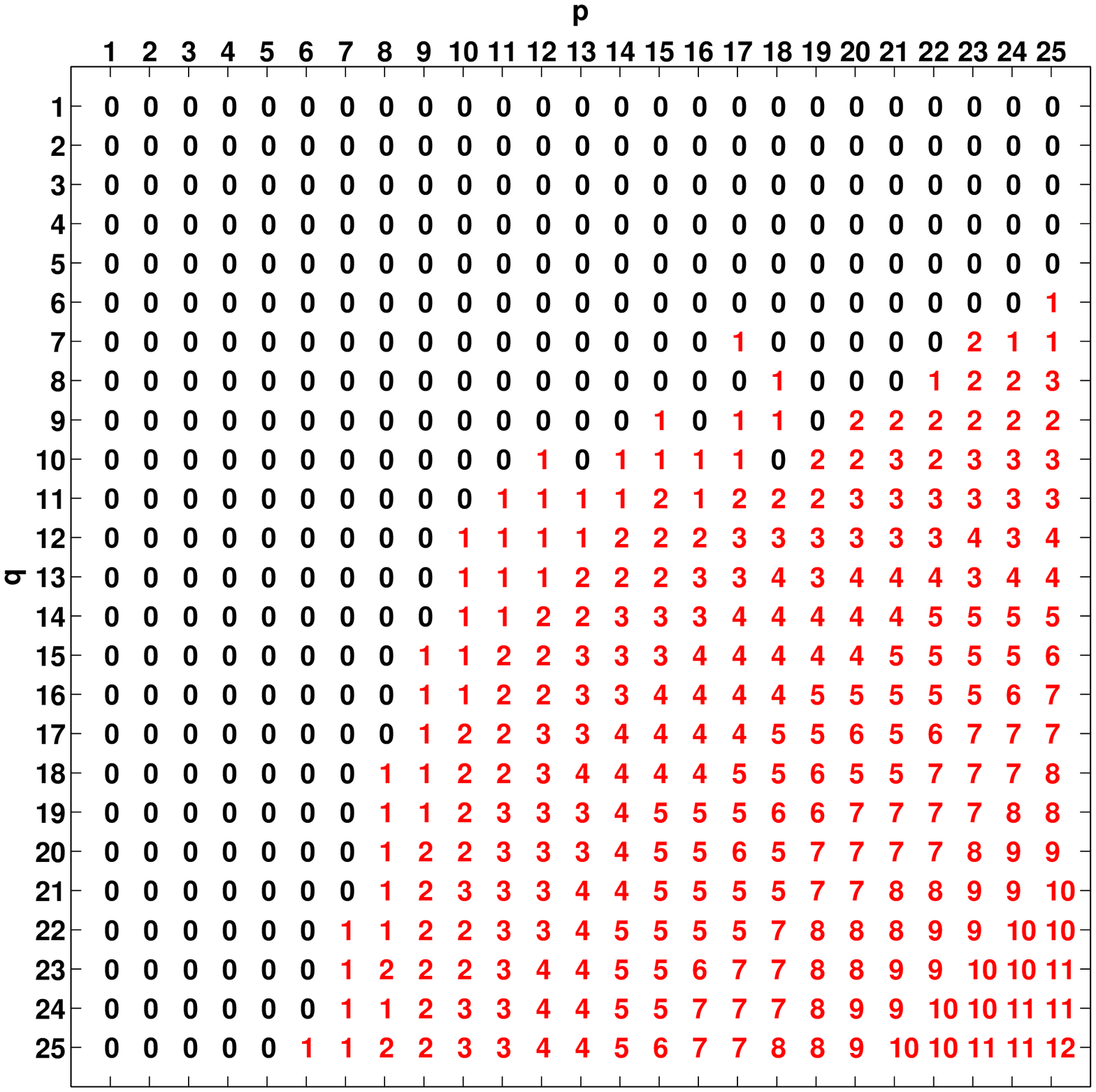}
\caption{Froissart table of the Tau solution $y_{150}$ with $\text{\rm tol}=10^{-5}$, $1\leq p,q\leq 25$.}
 \label{FroissartTable1}
\end{center}
\end{figure}
\begin{figure}[hbt]
\begin{center}
\includegraphics[width=1\textwidth]{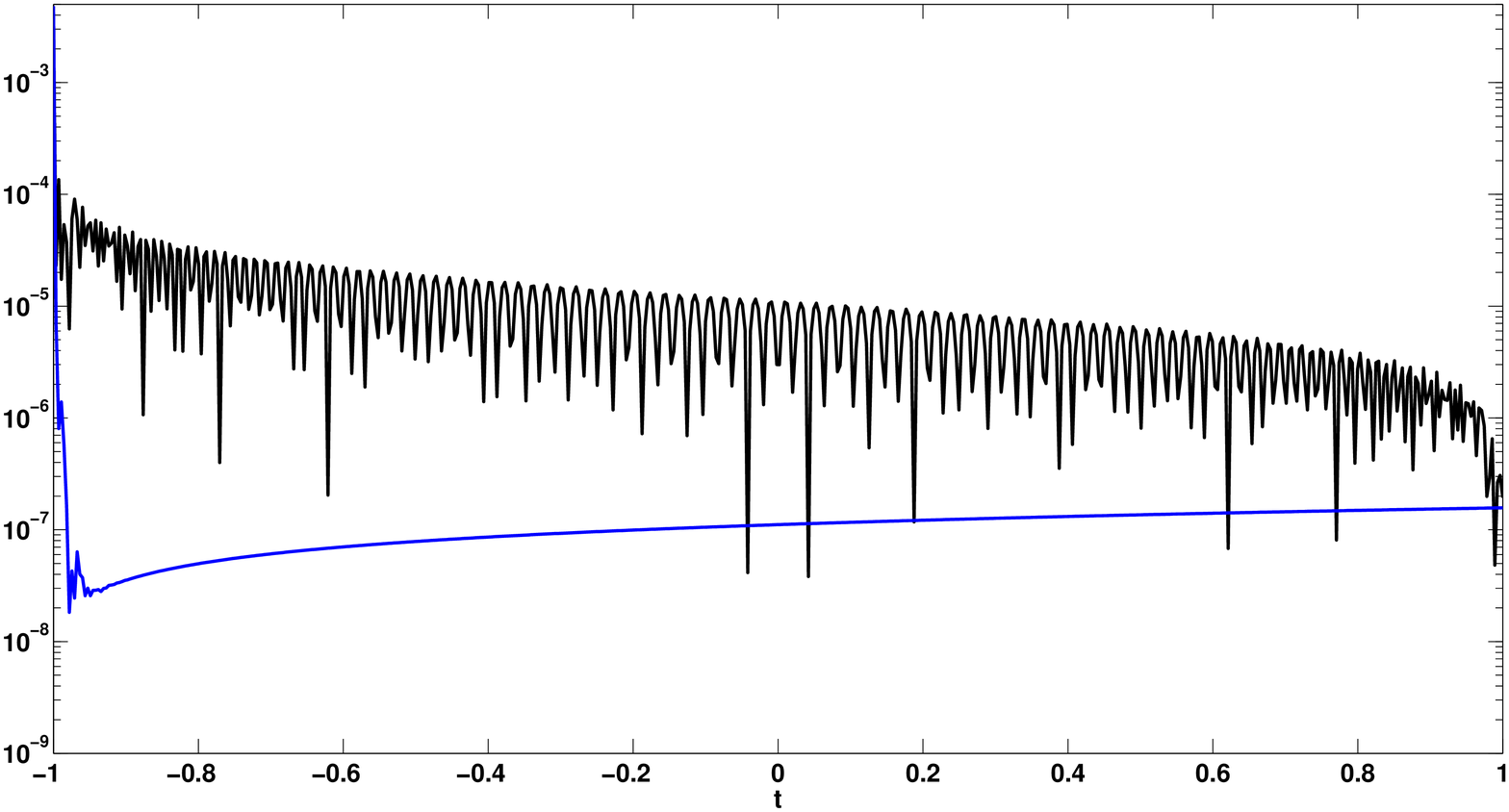}
\caption{Absolute error in Tau solution $y_{150}$, in black. Absolute error in the filter $\Phi_{10,10}(y_{150})$, in blue.}
 \label{fig:ErroDeFiltragem1}
\end{center}
\end{figure}

We present the Froissart table of the Tau solution $y_{150}$ with $\text{\rm tol}=10^{-5}$ and $1\leq p,q\leq 25$, in Figure \ref{FroissartTable1}.
In order to find a ``good'' filter, we look for one filter in the white region of the Froissart table. By other words we look for an $\Phi_{p,q}(y_{150};t)$
such that $n_{p,q}=0$ and that improves the Tau approximation. For example, if we look for a good diagonal filter $\Phi_{p,p}(y_{150};t)$, $p=1,2,\ldots,25$,
by inspecting the Froissart table, we can see that $\Phi_{10,10}(y_{150};t)$ is the last filter in the white region. All filters $\Phi_{p,p}(y_{150};t)$,
$p\geq 10$ have Froissart doublets and almost all have Froissart doublets on the real segment $[-1,1]$.

In Figure \ref{fig:ErroDeFiltragem1}, we illustrate ours results filtering the Tau solution $y_{150}$ of the example $1$ using the Pad\'e approximation $\Phi_{10,10}(y_{150};t)$. The absolute error $|y(t)-y_{150}(t)|$ is represented with a black line while the absolute error $|y(t)-\Phi_{10,10}(y_{150};t)|$ of the filter is
represented with a blue line.


\subsection{Estimation of singularities}
Another application of this filtering process is the estimation of singularities of the solution of a given differential problem. In fact, the poles of certain sequences of Pad\'e approximants from orthogonal series allow to estimate the singularities \cite{buslaev06,buslaev09}. However, we remark again, that our Pad\'e approximants are computed with the Tau coefficients, not with the orthogonal expansion coefficients.

For particular cases of orthogonal polynomial bases we can build particular formulas for poles $\lambda_{p}$, of the sequence of Pad\'e approximants $\Phi_{p,1},\ p=0,1,\ldots$ and for poles $\lambda_{p^{\pm}}$, of the sequence $\Phi_{p,2},\ p=0,1,\ldots$. The next results are consequences of corollary \ref{ChebyPadeApprox} and  of corollary \ref{LegPadeApprox}.

\begin{corollary}
Let $y(t)=\sum_{k=0}^\infty c_k T_k(t)$ be a formal Chebyshev series and $\Phi_{p,q}$
be its normalized Chebyshev-Pad\'e approximant of type $(p,q)$, then
\begin{itemize}
	\item[(a)] $\forall p\geq 1$ such that $c_{p+1}\neq 0$, the approximant $\Phi_{p,1}(y;t)$ has the pole 
\begin{equation}\label{polos1}
\lambda_{p}=\frac{c'_{p}+c_{p+2}}{2c_{p+1}}	
\end{equation}
	\item[(b)] $\forall p\geq 1$ such that the determinant
\[ \underline{\Delta} = \left|
\begin{array}{cc}
	c_{p+1} & c'_{p}+c_{p+2} \\
	c_{p+2} & c_{p+1}+c_{p+3}	
\end{array} \right| \neq 0 \]	
the approximant $\Phi_{p,2}$ has the poles
\begin{equation}\label{polos2}
\lambda_{p^\pm}=\frac{-b_{1}\pm \sqrt{b_{1}^{2}-8(b_{0}-1)}}{4}
\end{equation}
\noindent where $c'_p,\ b_{0}$ and $b_{1}$ are given in corollary \ref{ChebyPadeApprox} (b).
\end{itemize}

\end{corollary}

\begin{corollary}
Let $y(t)=\sum_{k=0}^\infty c_k P_k(t)$ be a formal Legendre series and $\Phi_{p,q}$
be its normalized Legendre-Pad\'e approximant of type $(p,q)$, then
\begin{itemize}
	\item[(a)] $\forall p\geq 1$ such that $c_{p+1}\neq 0$, the approximant $\Phi_{p,1}(y;t)$ has the pole 
\begin{equation}\label{Legpolos1}
\lambda_{p}=\frac{c_{p}+c_{p+2}}{2c_{p+1}}	
\end{equation}
	\item[(b)] $\forall p\geq 1$ such that the determinant
\[ \underline{\Delta} = \left|
\begin{array}{cc}
	c_{p+1} & h_{p+1,1} \\
	c_{p+2} & h_{p+2,1}	
\end{array} \right| \neq 0 \]	
\noindent the approximant $\Phi_{p,2}$ has the poles
\begin{equation}\label{Legpolos2}
\lambda_{p^\pm}=\frac{-b_{1}\pm \sqrt{b_{1}^{2}-3(b_{0}-1)}}{3}
\end{equation}
\noindent where $h_{k,1},\ b_{0}$ and $b_{1}$ are given in corollary \ref{LegPadeApprox} (b).
\end{itemize}

\end{corollary}

Thus, it is possible to compute the poles of the Pad\'e sequences $\Phi_{p,1}^{(n)}$ and $\Phi_{p,2}^{(n)}$, $p=0,1,\ldots$ only using the Tau coefficients on relations (\ref{polos1}-\ref{Legpolos2}) without computing Pad\'e approximants. 

Applying the relation \eqref{polos1} to the Tau solution of the Example \ref{exampleSQRT}, $y^{(150)}$, we obtained the results shown in Table \ref{tab:poles1}. We can see that all poles lie on the branch cut of $y$ and come closer to the singularity $\zeta =-1$, with increasing $p$. We remark that, in that case, since the exact Fourier coefficients $c_p=2(-1)^{p+1}/(4p^2-1)$ are known, we have access to the exact formula $\lambda_p=-1-3/(p-1/2)(p+5/2)$, easily obtained by substitution on \eqref{polos1}.

\begin{table}[htb]
  \centering

    \begin{tabular}{| c || c | c | c | c | c | }\hline
		$p$& $5$ & $45$ & $85$ & $105$ &  $145$  \\ \hline \hline
$\lambda_{p}$	& $-1.088889$ &  $-1.001419$ &  $-1.000406$ & $-1.000267$ &  $-1.000141$  \\ 
    \hline
    \end{tabular}
		  \caption{Poles of $\Phi_{p,1}^{(150)}$ for somes values of $p$.}
  \label{tab:poles1}
\end{table}

\section{Algebraic properties of filters}

In previous sections we saw how to obtain a rational approximation of a series, even when its coefficients are unknown. In this section we present some numerical properties of Frobenius-Pad\'e filters, justifying the observed results in numerical experiments.  

Our first property results from the following corollary of proposition \ref{Prop:ReqRel}:

\begin{corollary}\label{corol:filter}
Let 
\[ 
y = \sum_{i=0}^{\infty} c_i \phi_i \quad \text{and}\quad z = \sum_{i=0}^{\infty} d_i \phi_i
\] 
be two formal series where $\left\{\phi_i\right\}_{i\geq 0}$ is an orthogonal polynomial base satisfying \eqref{RecRel}. Let $H_y^{[p/q]}$ and $H_z^{[p/q]}$ be the matrices defined in \eqref{AFPsist1}, with subscript $y$ and $z$ identifying the respective series and with similar notation for other matrices and vectors in \eqref{AFPsist1} and \eqref{AFPsist2}. 

If 
\[ d_i=\rho\, c_i,\ i=0,1,\ldots,n \]
with a constant $0< |\rho|< \infty$ and if $H_y^{[p/q]}$ is regular, then 
$\forall p,q\in \mathbb{N}$ such that $p+2q\leq n$
we have:
\begin{itemize}
	\item[(a)] $H_z^{[p/q]}=\rho\, H_y^{[p/q]},\quad G_z^{[p/q]}=\rho\, G_y^{[p/q]},\quad  h_z^{[p/q]}=\rho\, h_y^{[p/q]}\ $ \\ and $\ g_z^{[p/q]}=\rho\, g_y^{[p/q]}$;
	\item[(b)] $H_z^{[p/q]}$ is regular;
	\item[(c)] $N_{p,q}(z)=N_{p,q}(y)$ and $D_{p,q}(z)=\rho\, D_{p,q}(y)$, where $N$ and $D$ are the numerator and denominator polynomials introduced in \eqref{phi};
	\item[(d)] $\Phi_{p,q}(z)=\rho\, \Phi_{p,q}(y)$ 
\end{itemize}
 
\end{corollary}

The proof of (a) follows by induction over $j$ in \eqref{RecRel} and then (b) and (c) results from (a), and (d) results from (b) and (c).  
 
From this results we get the following corollary, related to the Frobenius-Pad\'{e} filtering of spectral approximation of Fourier series, which explains the good behaviour of the filtering process. 

\begin{corollary}\label{corol:poles}
Let $p,q,m,n\in\mathbb{N}_0$ and $\Phi_{p,q}(y_n)$ be a $(p,q)$ Pad\'e filter from $y_n$. Suppose that for some finite and non null constant $\rho\in\mathbb{R}$, the coefficients $c_k^{(n)}=\rho\, c_k,\ k=0,\ldots,m$ and $m\leq n$.
\begin{itemize}
	\item[(a)] If $p+2q\leq m$ and $H_{y_n}^{[p/q]}$ is regular then $\Phi_{p,q}(y_n)=\rho\, \Phi_{p,q}(y)$;
	\item[(b)] If $p\leq m-2$ then $\lambda_{p}^{(n)}=\lambda_{p}$, that is, the pole of the filter $\Phi_{p,1}(y_n)$ coincides with the pole of $\Phi_{p,1}(y)$;
	\item[(c)] If $p\leq m-4$ then $\lambda_{p^\pm}^{(n)}=\lambda_{p^\pm}$, that is, the poles of the filter $\Phi_{p,2}(y_n)$ coincide with the poles of $\Phi_{p,2}(y)$ 
\end{itemize}  
\end{corollary}

So, if we have a set of numerical approximations $c_k^{(n)}\approx c_k,\ k=0\ldots,m$ and if all of those approximations have the same relative errors $\delta=(c_k-c_k^{(n)})/c_k,\ k=0\ldots,m$, then corollaries \ref{corol:filter} and \ref{corol:poles} old with $\rho=1-\delta$ and our filter process, working with $c_k^{(n)},\ k=0\ldots,m$, will give rational approximants with relative error $\delta$ and with the same poles as if we work with the exact coefficients $c_k,\ k=0\ldots,m$.

\section{Numerical example}

In the next example, we will test this filtering procedure to  Legendre-Tau solutions of a differential equation with boundary conditions. Furthermore, since the differential equation depends on a parameter, that allows to control
the rate of convergence of the Legendre-Tau method, we can test the behavior of the filters when applied to problems with different rates of convergence.

\begin{examp}\label{genfunLeg}
Let us consider the family of  functions
\[ y(t)=\frac{1-\alpha^2}{(1+\alpha^2-2\alpha t)^{3/2}},\ t\ ]-1, 1[ \]
depending on the real parameter $\alpha$ and whose Legendre series representation
\[ y(t)=\sum_{k=0}^{\infty} (2k+1)\alpha^k P_k(t),\ t\in [-1, 1] \]
can be derived from the generating function of Legendre polynomials \cite{Abramowitz65}. For $\alpha\neq 0$, $y$ has  branch points at $\zeta=\frac{1}{2}(\alpha+\frac{1}{\alpha})$ and at $\infty$, and furthermore, $\zeta$ is the closest singularity of the interval of orthogonality $[-1,1]$. 

For our purpose, we can define $y$ as the solution of the boundary value pro\-blem
\begin{equation} \label{examp2}
\left\{
	\begin{array}{ll}
	(1+\alpha^2-2\alpha t)^{2}y''(t)-15\alpha^2 y(t)=0, & t\in [-1, 1] \\
	y(-1)=\frac{1-\alpha}{(1+\alpha)^2},\ y(1)=\frac{1+\alpha}{(1-\alpha)^2} 
	\end{array}
\right.
\end{equation} 

\noindent The matricial form of the operator $D$, introduced in (\ref{piphi}), associated to this differential problem  is given by:
\[ \Pi_{\phi} = \eta_{\phi}^{2}((1+\alpha^2)I-2\alpha \mu_{\phi})^2-15\alpha^2 I \]
where $I$ is the infinity identity matrix and with $\phi$ being the Legendre polynomials.

\begin{figure}[hbt]
\begin{center}
\includegraphics[width=1\textwidth]{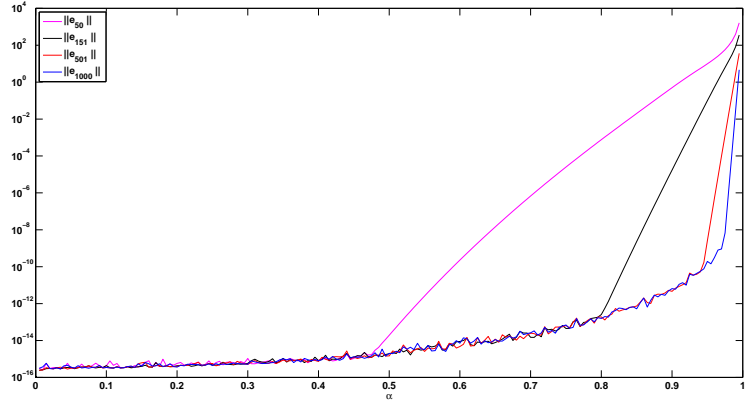}
\caption{$w$-norm of the functions errors $e_{n}$, $n=50, \ 151, \  501$ and $1000$, for values of $\alpha\in ]0,1[$.}
 \label{Errotauleg1a}
\end{center}
\end{figure}

The rate of convergence of the Tau method applied to problem (\ref{examp2}) depends on the parameter $\alpha$, exhibiting slow rate of convergence for values of $\alpha$ nearby $1$. To illustrate this behavior we show
in Figure \ref{Errotauleg1a}, the weighted-norm $||e_{n}||_{w}$ of the errors of four Legendre-Tau solutions $y_{n}$, $n=50,151,501$ and $1000$, for values $\alpha\in]0,1[$.  We can see that for $\alpha<0.5$ it is  enough to compute $y_{50}$ to get an error of order of the machine precision. For values of $\alpha>0.5$, we need to increase the order of the Legendre Tau solutions, to get a reasonable approximation and the machine precision is lost.

\begin{figure}[hbt]
\begin{center}
\includegraphics[width=1\textwidth]{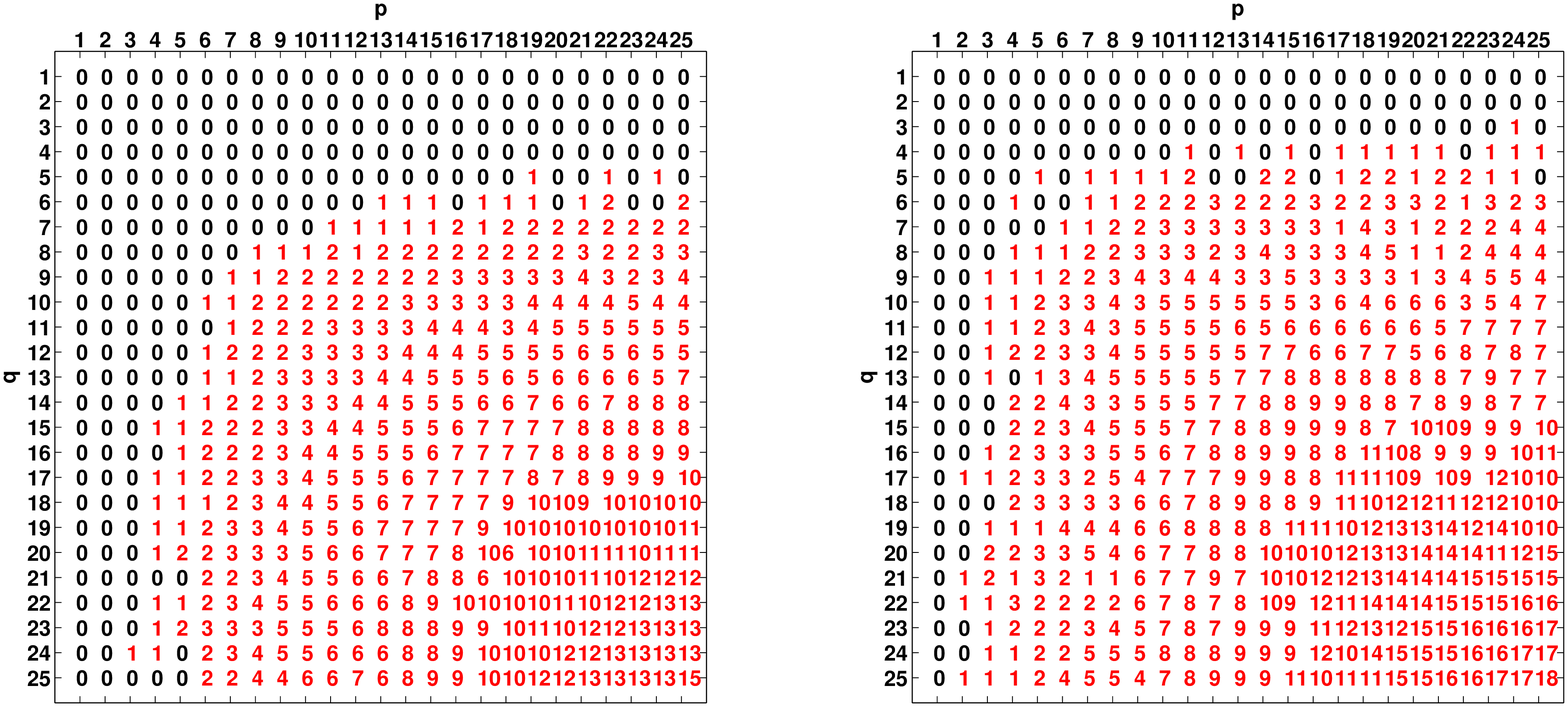}
\caption{Froissart Tables, computed with $tol=10^{-5}$, of $y_{150}$ with $\alpha=0.9$  in the left table and for $y_{1000}$ with $\alpha=0.99$ in the right.}
 \label{legendrefroissart1}
\end{center}
\end{figure}

\begin{figure}[hbt] 
\begin{center}
\includegraphics[width=1\textwidth]{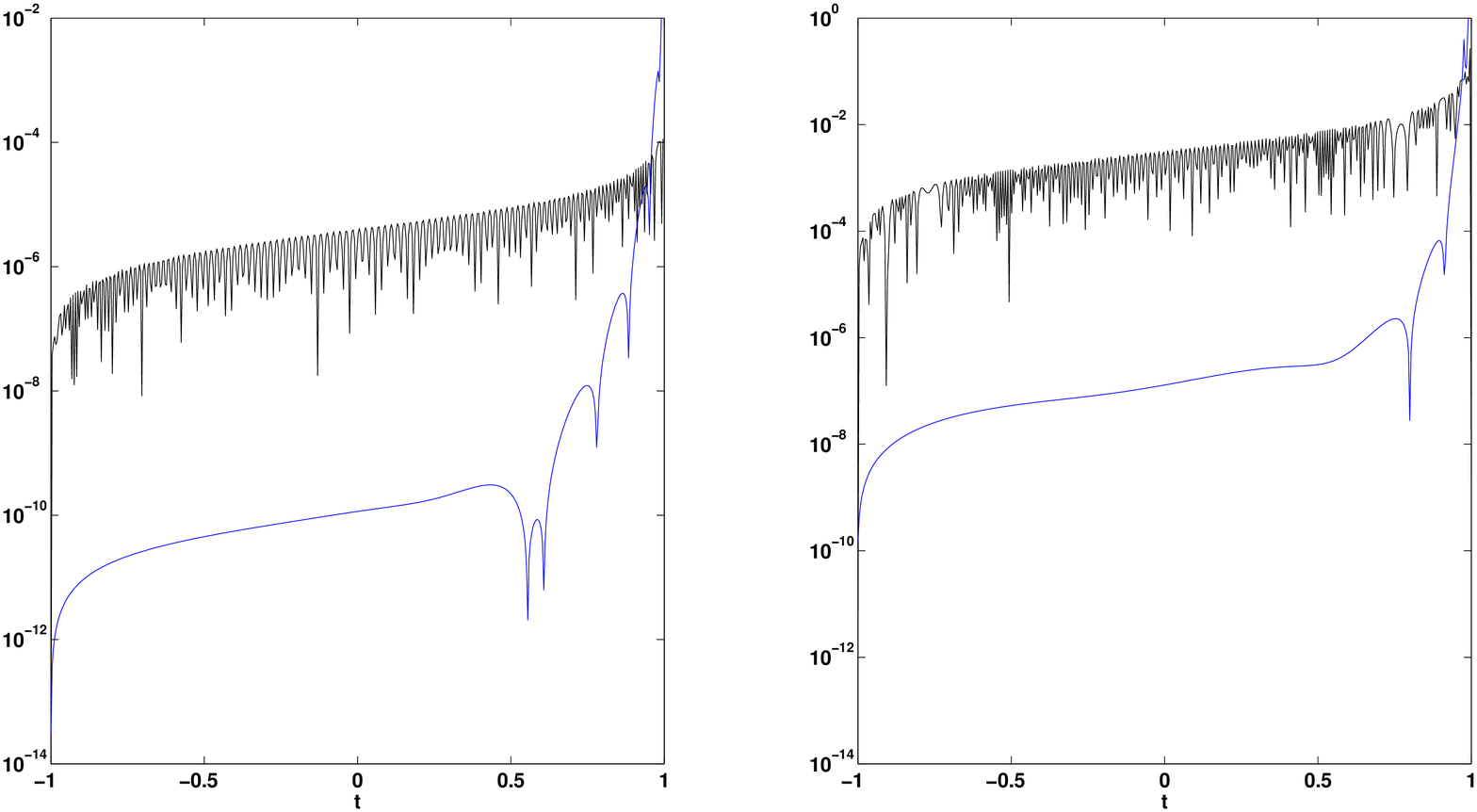}
\caption{ Left image: Absolute error of theTau solution, $y_{150}$, of the problem with $\alpha=0.9$  (black line) and  absolute error of the filter
$\Phi_{7,7}(y_{150})$ (blue line). Right image: Absolute error of the Tau solution, $y_{1000}$, of the problem with $\alpha=0.99$  (black line) and absolute error of the filter
$\Phi_{6,6}(y_{1000})$ (blue line).}
\label{errosfiltragem}
\end{center}
\end{figure}

In order to test our filtering method, we computed the Legendre Tau solution $y_{150}$ of (\ref{examp2}) with $\alpha=0.9$ and
the Legendre Tau solution $y_{1000}$ of the same problem with $\alpha=0.99$. Proceeding in analogous way to the example
\ref{exampleSQRT}, we took for a ``good'' filter the diagonal Legendre Pad\'e approximant $\Phi_{p,p}(y_{n};t)$  in the
white region of the Froissart Table  for which  $p$ is maximum. The Figure \ref{legendrefroissart1} shows the Froissart tables (with $tol=10^{-5}$) of
$y_{150}$, $\alpha=0.9$ (left table) and of $y_{1000}$, $\alpha=0.99$ (right table). Inspecting the tables we see that $\Phi_{7,7}(y_{150};t)$
is a ``good'' diagonal filter for the first problem while for the second problem we must choose $\Phi_{6,6}(y_{1000};t)$. In Figure \ref{errosfiltragem} we show the absolute errors of the tau solutions $y_{150}$ and  $y_{1000}$ and the absolute errors of their filters, $\Phi_{7,7}(y_{150})$
and $\Phi_{6,6}(y_{1000})$, respectively. In both cases, the filters improve the Legendre Tau approximations for values of $t\in[-1,1]$ that are not close of $t=1$.

In order to estimate $\zeta$ we can use the relation \eqref{Legpolos1} to compute the zeros
of the filters $\Phi_{p,1}(y_{n})$ and use them as approximants of $\zeta$. However, this problem it is a differential equation with  boundary conditions and we did not get a relation between the Tau coefficients and the Fourier coefficients, as in example \ref{exampleSQRT}. In fact, the relative error of the Tau coefficients are not constant and we need proceed carefully, because the poles of $\Phi_{p,1}(y_{n})$ have not the same behavior of the poles of $\Phi_{p,1}(y)$. %
\end{examp}

\section{Conclusions}

Our numerical experiments reveal that it is possible to improve the Tau solutions approximations using Pad\'e approximation. The noise introduced on the Tau coefficients in the numerical computation of Pad\'e approximants yields the occurrence of Froissart doublets for high order rational approximants. The Froissart table, introduced in this work, reveals to be an efficient tool to find a good filter of the Tau solution.

This filtering method also allows to estimate singularities of exact solutions, since the computation of the poles of $\Phi_{p,1}(y_{n})$ and $\Phi_{p,2}(y_{n})$ can be computed  using only the Tau coefficients.

Some algebraic properties of the filtering process were introduced, justifying the good properties of the filtered solutions.


\bibliography{bibliografia}

\end{document}